\documentclass[12pt]{paper}

\usepackage{amssymb,amsfonts,amsthm,amsmath}

\makeatletter
\newcommand{\subjclass}[2][1991]{%
  \let\@oldtitle\@title%
  \gdef\@title{\@oldtitle\footnotetext{#1 \emph{Mathematics subject classification.} #2}}%
}
\makeatother

\usepackage{graphicx, verbatim}
\usepackage{xcolor}

\usepackage[a4paper, hmargin=3cm, vmargin={4cm, 4cm}]{geometry}
\usepackage{fancyhdr}

\newtheorem{theorem}{Theorem}
\newtheorem{lemma}[theorem]{Lemma}

\newtheorem{proposition}[theorem]{Proposition}

\newtheorem{claim}[theorem]{Claim}

\newcommand{\bR}{\mathbb R}
\newcommand{\bC}{\mathbb C}
\newcommand{\bZ}{\mathbb Z}
\newcommand{\bD}{\mathbb D}

\newcommand{\bE}{\mathbb E}

\newcommand{\eps}{\varepsilon}
\newcommand{\la}{\lambda}
\renewcommand{\phi}{\varphi}

\begin{document}

\title{Taylor coefficients and zeroes \\ of entire functions of exponential type
}

\author{Lior Hadassi, Mikhail Sodin}

\maketitle
\hfill{\sf\today}

\begin{abstract}
Let $F$ be an entire function of exponential type represented by the Taylor series
\[
F(z) = \sum_{n\ge 0} \omega_n \frac{z^n}{n!}
\]
with unimodular coefficients $|\omega_n|=1$. We show that {\em either the counting function $n_F(r)$
of zeroes of $F$ grows linearly at infinity,
or $F$ is an exponential function}.
The same conclusion holds if only a positive asymptotic proportion of the coefficients
$\omega_n$ is unimodular. This significantly extends a classical result of Carlson (1915).

The second result requires less from the coefficient sequence $\omega$, but more from the counting function of zeroes $n_F$. Assuming that $0<c\le |\omega_n| \le C <\infty$, $n\in\mathbb Z_+$, we show that {\em 
$n_F(r) = o(\sqrt{r})$ as $r\to\infty$, implies that $F$ is an exponential function}. The same conclusion holds if, for some $\alpha<1/2$, $n_F(r_j)=O(r_j^{\alpha})$ only along a sequence $r_j\to\infty$.
Furthermore, this conclusion ceases to hold if $n_F(r)=O(\sqrt r)$ as $r\to\infty$.
\end{abstract}

\section{Introduction}
In this work we study the following question: how a sequence of multipliers $\omega\colon \bZ_+~\to~\bC$
affects the number of zeroes of the entire function
\[ F (z) = \sum_{n\ge 0} \omega_n\, \frac{z^n}{n!}\,? \]
We will assume that $F$ is an entire function of exponential type (EFET, for short), 
that is, $\log M_{F}(R) = O(R)$ as $R\to\infty$, where
$M_{F}(R)=\displaystyle \max_{|z|=R} |F(z)|$, and therefore, the number of zeroes of
$F$ in the disk $\{|z|\le R\}$ grows at most linearly.

The relationship between the properties of the sequence $(\omega_n)$ and the zero distribution of 
$F$ has been studied by many authors (Littlewood and Offord, Levin, Nassif, and others)
under the assumption that the sequence $(\omega_n)$ possesses a special structure. The most typical 
instances are random i.i.d. sequences and, more generally, random stationary sequences, sequences of 
arithmetic nature, and almost-periodic sequences. We refer the reader to the recent works by
Borichev, Nishry, and Sodin~\cite{BNS} and Benatar, Borichev, and Sodin~\cite{BBS1, BBS2}, which also contain 
the references to what was done previously.

Here, our departing point is a classical result of Carlson~\cite{Carlson}, who showed that {\em if 
$|\omega_n|=1$ for $n\in\bZ_+$, and $F$ has few zeroes in the sense that
\begin{equation}\label{eq:Carlson}
\sum \frac1{|\la|} < \infty\,,
\end{equation}
where the summation is taken over the zero set of $F$,
then $F$ is an exponential function, that is,  $F(z)=\Theta e^{\alpha z}$ with $|\Theta|=|\alpha|=1$}.
This is the only result in this direction that we are aware of, which does not assume a special structure for
$\omega_n$.

\medskip In this work, we extend Carlson's result in two directions.
Our first result shows that~\eqref{eq:Carlson} can be replaced with a much weaker assumption.
Let $n_F(R)$ be a number of zeroes of $F$ in the disk $\{|z|\le R\}$,
counted with multiplicities.
\begin{theorem}\label{thm:main1}
Let
\[
F(z) = \sum_{n\ge 0} \omega_n \frac{z^n}{n!}
\]
with $|\omega_n|=1$, $n\in\bZ_+$. Then there exists a positive numerical constant $c$ such that
either $F$ is an exponential function, or
\begin{equation}\label{eq:local}
\liminf_{R\to\infty} \frac{n_F(R)}R \ge c>0.
\end{equation}
\end{theorem}
Note that the Carlson assumption~\eqref{eq:Carlson} corresponds to convergence of the integral
\[
\int^\infty \frac{n_F(t)}t\, {\rm d}t <\infty.
\]
It is also worth mentioning that we do not know a sharp value of the constant $c$ in Theorem~\ref{thm:main1}.

\medskip
Actually, we will prove more, requiring only that the set $\{n\colon |\omega_n|=1\}$ has a positive lower density, and that
$F$ has a finite exponential type, see Theorem~\ref{thm:main} below. Note that in the end of his note Carlson also mentions that his result holds when $F$ has a finite exponential type and the  set $\{n\colon |\omega_n|=1\}$ has a positive lower density.

\medskip
The proof of Theorem~\ref{thm:main1} uses several novel tools comparing with the ones used by Carlson. Let us briefly
recall his argument. Condition~\eqref{eq:Carlson} yields that $F(z)=e^{\alpha z} G(z)$, where $|\alpha|=1$ and $G$ is an entire function of zero exponential type.  Rotating the coordinates, assume that $\alpha =1$, i.e., $F(z)=e^z G(z)$.
Let $\widehat F$ be the Laplace transform of $F$. It is analytic in $\overline{\bC}\setminus\{1\}$, vanishes at infinity,
and has the expansion
\[
\widehat{F}(w) = \sum_{n\ge 0} \frac{\omega_n}{w^{n+1}}\,.
\]
Then, its coefficients $\omega_n$ can be interpolated by an entire function $g$ of zero exponential type, that is,
$g(n)=\omega_n$, $n\in\bZ_+$. Then the entire function
$\widetilde{g}(\la)=1-g(\la)\overline{g}(\overline \la)$ has zero exponential type and vanishes on $\bZ_+$.
Thus, by a classical uniqueness theorem (which is also due to Carlson), $\widetilde g$ is the zero function, i.e.,
$g$ is a unimodular constant, that is, $\omega$ is a constant sequence, and $F(z)=Ce^z$, $|C|=1$.

Note that Carlson's theorem was extended by Eremenko and Ostrovskii~\cite[Theorem~1]{EO}
in a direction different from the one pursued here. They showed that  if $|\omega_n|=1$ for $n\in\bZ_+$, then $F$ cannot
decay exponentially fast on any ray $\{\arg(z)=\theta\}$. Clearly, this yields Carlson's theorem, since the representation
$F(z)=e^{\alpha z} G(z)$, where $G$ has zero exponential type and $|\alpha|=1$, immediately yields that $F$ decays exponentially fast in an open half-plane. Their argument also was based on the consideration of the Laplace transform
$\widehat{F}$.

The main difficulty in extending Carlson's idea to our case is that the negation of our assumption~\eqref{eq:local} 
yields only local estimates on the function $G$. These local estimates do not  allow us to use the Laplace transform and classical results on the interpolating function. To circumvent this obstacle, instead of
interpolating the sequence $\omega_n$, we interpolate the sequence $\omega_n \bar\omega_{n+h}$, where
$h\in\bZ_+$ is a parameter, which needs to be adjusted. It turns out that the local information on $F$ that Theorem~\ref{thm:main1} provides, is sufficient to efficiently bound the interpolating function. In addition,
we will need a small piece of additive combinatorics to deal with the case when  the set $\{n\colon |\omega_n|=1\}$ only
has a positive lower density.

It is worth mentioning that the techniques developed for the proof of Theorem~\ref{thm:main1} allow one 
also to treat other closely related families of entire functions, such as
\[
\sum_{n\ge 0} \omega_n\, \frac{z^n}{\sqrt {n!}} \quad  \text{and} \quad  \sum_{n\ge 0} \omega_n\, \frac{z^n}{\Gamma (n+\alpha)}\,.
\]

\medskip
Our second result requires less from the coefficient sequence $\omega$, but more from the counting function of zeroes $n_F$.
\begin{theorem}\label{thm:main2}
Let
\[
F(z) = \sum_{n\ge 0} \omega_n \frac{z^n}{n!}
\]
with $ 0<c \le |\omega_n| \le C < \infty$, $n\in\bZ_+$. Suppose that 

\smallskip\noindent {\rm (i)\ } either
\begin{equation}\label{eq:n_F}
\limsup_{R\to\infty}\, \frac{n_F(R)}{\sqrt R} = 0,
\end{equation}

\smallskip\noindent {\rm (ii)\ }
or 
\begin{equation}\label{eq:n_F-alpha}
\liminf_{R\to\infty}\, \frac{n_F(R)}{R^{\frac12 - \eps}} = 0
\end{equation}
with some $\eps>0$.

\noindent Then $F$ is an exponential function.
\end{theorem}

The proof of both parts starts with the observation that, by the boundedness from below and from above of
the sequence $(\omega_n)$, we have
\[
R - \tfrac12\, \log R - O(1) < \log M_F(R) < R + O(1), 
\qquad R\to \infty.
\]
The rest of the argument is purely subharmonic.
The proof of part (i) is quite simple. In fact, it yields 
a somewhat stronger result (see Remark at the end of Section~3.1).
The proof of part (ii) (suggested by Fedor Nazarov)
is more intricate, and we do not know whether one can replace~\eqref{eq:n_F-alpha}
by 
\begin{equation}\label{eq:n_F-1/2}
\liminf_{R\to\infty}\, \frac{n_F(R)}{\sqrt{R}} = 0.
\end{equation}
 
We conclude the paper with two examples. 
First, we will show that the entire function
\[
F(z) = \sum_{n\ge 0}\, (\cos\sqrt{n}+2)\,\frac{z^n}{n!}
\]
can be written in the form $F(z)=e^z G(z)$, where $G$ is an entire function of order $1/2$ and mean type. Therefore,
$n_F(R)=O(\sqrt R)$ as $R\to\infty$; i.e., the exponent $\tfrac12$ in part (i) of Theorem~\ref{thm:main2} cannot
be improved. The second example will illustrate why the method used in the proof of Theorem~\ref{thm:main2}(ii) 
ceases to work under a weaker assumption~\eqref{eq:n_F-1/2}.

\paragraph{Acknowledgement} \mbox{}

\smallskip\noindent 
Fedor Nazarov kindly helped us with a proof of part (ii) of Theorem~\ref{thm:main2} and provided
the aforementioned example illustrating this part, but unfortunately declined to co-author this paper.
Lev Buhovskii, Alex Eremenko, and Alon Nishry read most of this paper, and provided us with many
comments and remarks, which we took into account.  We are sincerely grateful to them all.

The work of both authors was supported by the ISF Grant 1288/21. 

\paragraph{Notation}$\mbox{}$

\begin{itemize}
\item $M_F(R)=\max_{|z|=R}|F(z)|$.
\item $n_F(R)$ denotes the number of zeroes of
$F$ in the disk $\{|z|\le R\}$ counted with multiplicities.
\item {\sf EFET} means entire function of exponential type, i.e.,
$|F(z)|\le Ae^{B|z|}$, $z\in\mathbb C$. \newline
By ${\sf EFET}(\sigma)$ we denote the class of entire functions of exponential type at most $\sigma$,
that is, for each $\delta>0$, $|F(z)|\le A_\delta e^{(\sigma+\delta)|z|}$, $z\in\mathbb C$. 
\item $f^*(z)=\overline{f(\bar z)}$.
\item $C$ and $c$ (with or without indices) are positive constants.
One can think that the constant $C$ is large (in particular, $C \ge 1$), while the constant
$c$ is small (in particular, $c \le 1$).
The values of these constant may change from line
to line. Within each section, we start a new sequence of indices for these constants.

\item $A\lesssim B$ means $A\le C\cdot B$, while $A\gtrsim B$ means $A\ge cB$.
We use the standard $O$-notation and the notation $\lesssim $ with interchangeable meaning.

\item $\ll$ stands for  “sufficiently smaller than” and means that $A<cB$ with sufficiently small positive $c$,
while $\gg$ stands for  “sufficiently larger than” and means  that $A > C B$ with a very large positive $C$.

\item For a set of non-negative integers $X\subset \mathbb Z_+$, $d(X)$ denotes its lower
density:
\[
d(X) = \liminf_{N\to\infty}\, \frac{|X\cap [1, N]|}N\,.
\]

\item For a set $Y\subset\bC$, $Y_{+t}=\{z\colon \text{dist}(z, Y)<t\}$ denotes the $t$-neighbourhood
of $Y$.

\item As usual, $\log_+ x = \max\{ \log x, 0\}$ and $\log_- x = \max\{ -\log x, 0\}$.

\end{itemize}

\section{Unimodular coefficients}
In this section we prove

\begin{theorem}\label{thm:main}
Let
\[
F(z) = \sum_{n\ge 0} \omega_n \frac{z^n}{n!}
\]
be an {\sf EFET}$(\sigma)$ such that the set $\Lambda_F=\{n\in\bZ_+\colon |\omega_n|=1\}$ has a
positive lower density $d>0$. Then there exists $\chi$ depending only on the lower density $d$ and type $\sigma$
such that either
\[
\liminf_{R\to\infty} \frac{n_F(R)}R \ge \chi > 0
\]
or $F(z)= \Theta\, e^{\alpha z}$ with $|\alpha|=|\Theta|=1$.
\end{theorem}

To avoid a large number of parameters in the proof of Theorem~\ref{thm:main}, 
we will employ the following convention:
\begin{itemize}
\item notation $o_\eps (1)$ will denote any positive quantity which depends only on a small positive parameter $\eps$
and tends to zero as $\eps\downarrow 0$, while $\omega_\eps(1)$ will denote any positive quantity which 
depends only on $\eps$ and tends to $+\infty$ as $\eps\downarrow 0$. 
\end{itemize}

\subsection{The self-correlation transform}

Given a Taylor series
\[
f(z) = \sum_{n\ge 0} a_n z^n
\]
and $h\in\bZ_+$, we set
\begin{align*}
f_h(z) &= \sum_{n\ge 0} a_n \bar a_{n+h} z^{2n} \\
&= \oint_{|s|=1} f(sz) f^*(\bar s z)
\Bigl( \frac{s}z \Bigr)^h\, \frac{{\rm d} s}{2\pi {\rm i}s}\,,
\end{align*}
where the integral is taken over the unit circle traversed counterclockwise.
Applying this transform to our function $F$, we get
\[
F_h(z) =
\sum_{n\ge 0} \omega_n \bar\omega_{n+h}\, \frac{z^{2n}}{n!(n+h)!}\,.
\]

\begin{lemma}\label{lemma-SCT}
Suppose that there exist $R \ge 1$ and $\alpha\in\mathbb C$ such that
\[
F(z) = e^{\alpha z} G(z)\,,
\]
\[
\log M_G(R) < o_\eps (1) R\,,
\]
and
\[
|\alpha| < 1 +o_\eps(1)\,.
\]
Then, for any $r\in [1, R]$, we have
\[
|F_h(re^{{\rm i}\theta})| \le r^{-h} e^{2r|\cos \theta| + o_\eps (1) R}\,.
\]
\end{lemma}

\noindent{\em Proof}: We have
\begin{align*}
|F_h(re^{{\rm i}\theta})| &= \Bigl|\, \oint_{|s|=1}
F(sre^{{\rm i}\theta}) \overline{ F(sre^{-{\rm i}\theta})}
\Bigl( \frac{s}{re^{{\rm i}\theta}} \Bigr)^h\, \frac{{\rm d} s}{2\pi {\rm i}s}  \,\Bigr| \\
&\le r^{-h}\, \max_{|s|=1} \bigl| e^{\alpha s r e^{{\rm i}\theta}} \cdot
e^{\bar \alpha \bar s r e^{{\rm i}\theta}}  \bigr|\, e^{o_\eps(1) R} \\
&=  r^{-h} e^{2r|\alpha|\,|\cos \theta| + o_\eps(1) R} 
\qquad (\text{since}\ \max_{|s|=1}\, | \alpha s + \bar\alpha \bar s | = 2|\alpha|) \\
&\le r^{-h} e^{2r|\cos \theta| + o_\eps(1) R},
\end{align*}
proving the lemma. \hfill $\Box$

\medskip
It will be more convenient to work with the function
\[
F_h^\sharp (z) = F_h(\sqrt z) = \sum_{n\ge 0} \omega_n \bar\omega_{n+h} \frac{z^n}{n!(n+h)!}\,.
\]
By Lemma~\ref{lemma-SCT}, for $z=re^{{\rm i}\theta}$, $1 \le r \le R^2$, $|\theta| \le \pi$,
we have
\[
|F_h^\sharp (re^{{\rm i}\theta})| \le r^{-h/2} e^{2\sqrt{r}\cos \frac{\theta}2
+ o_\eps(1) R}\,.
\]

\subsection{Analytic interpolation}

Given $0\le h \le R$, we set
\[
\kappa = \Bigl[\,  \frac{h}2\, \Bigr] + \frac12
\]
(as usual, $[x]$ denotes the integer part of $x$),
and denote by $\Omega_h$ the shifted semi-disk
\[
\Omega_h = \bigl\{
\lambda\colon \text{\ Re}(\lambda) >-\kappa, |\lambda+\kappa|<R-\kappa
\bigr\}.
\]
\begin{lemma}\label{lemma_interpol}
Let
\[
F(z)=\sum_{n\ge 0} \omega_n\, \frac{z^n}{n!}
\]
be an {\sf EFET}. Suppose that, for some $R>1$ and $\alpha\in\mathbb C$, we have
\[
F(z) = e^{\alpha z} G(z),
\]
\begin{equation}
\label{eq_G}  \log M_G(R) < o_\eps (1) R,
\end{equation}
\begin{equation}
\label{eq_alpha} |\alpha| < 1 + o_\eps(1).
\end{equation}
Then, given $h\in\bZ_+$, there exists a function $g_h$ holomorphic
in the half-plane $\{ \text{Re}(\lambda) > -(h+1)\}$, interpolating the sequence $\omega_n \bar\omega_{n+h}$, i.e.,
\[
g_h(n) = \omega_n \bar\omega_{n+h}\,, \qquad n\in\mathbb Z_+,
\]
and, for any $h\in [0, R]$, satisfying
\[
|g_h(\la)|\le e^{o_\eps(1)R +O(R^{-1}h^2)}\,, \qquad \la\in\bar\Omega_h,
\]
provided that $R=R_\eps$ is sufficiently big.
\end{lemma}

\noindent{\em Proof of Lemma~\ref{lemma_interpol}}:
Note that, for $n\in\bZ$,
\[
n! (n+h)!\, \oint_{|z|=\rho} \frac{F_h^\sharp (z)}{z^n}\, \frac{{\rm d}z}{2\pi {\rm i}z}
=
\begin{cases}
\omega_n \bar\omega_{n+h} & n\ge 0, \\
0 & n<0,
\end{cases}
\]
where the integral is taken over the circle of radius $\rho>0$ traversed counterclockwise.
Keeping this in mind, given $h\in\bZ_+$, we consider the function
\begin{equation}\label{eq:an_interpol}
g_h(\la) =  \Gamma (\la+1) \Gamma (\la+h+1)\, \int_{-{\rm i}\pi + 2\log r}^{{\rm i}\pi + 2\log r}
F_h^\sharp (e^w) e^{-\la w}\, \frac{{\rm d}w}{2\pi {\rm i}}\,,
\end{equation}
where the integration is taken over the vertical segment $\{w=2\log r + {\rm i}t\colon |t|\le\pi\}$,
and the parameter $r$ has to be chosen.
The integral on the RHS defines an entire function in $\la$ vanishing at the negative integers, therefore, the function $g_h$ is analytic in the right half-plane $\{{\rm Re}(\la)>-(h+1)\}$,
and interpolates $ \omega_n \bar\omega_{n+h} $ at $n\in\bZ_+$.

It remains to estimate $g_h$ in $\bar\Omega_h$, so in what follows we assume that $0\le h \le R$.
It will be convenient to shift the variable $\la$, letting
\[ s = \la + \kappa = \rho e^{{\rm i}\phi}, \ R_h = R-\kappa, \ 
D_+(R_h) = \{s\colon {\rm Re}(s)> 0, |s|< R_h \}\,. \]
In this notation, $\Omega_h = D_+(R_h)-\kappa$, and
\begin{align*}
\Gamma (\la +1) \Gamma (\la+h+1) &= \Gamma \bigl( s - \bigl[ \tfrac{h}2 \bigr] + \tfrac12 \bigr)
\Gamma \bigl( s + \bigl[ \tfrac{h}2 \bigr] + \tfrac12 \bigr) \\
&= \Gamma \bigl( s + \tfrac12 \bigr)^2\, \Pi_h(s)\,,
\end{align*}
where
\[
\Pi_h(s)=\prod_{1\le \ell \le [\frac{h}2]} \frac{s+(\ell - \frac12)}{s-(\ell-\frac12)}\,.
\]
At last, we set in~\eqref{eq:an_interpol} $r=R_h = R-\kappa$, and get
\begin{align}\label{eq:G_h}
g_h(\la) &= \Gamma \bigl( s + \tfrac12 \bigr)^2\, \Pi_h(s)
\int_{-{\rm i}\pi + 2\log R_h}^{{\rm i}\pi + 2\log R_h}
F_h^\sharp (e^w)e^{\kappa w} \cdot e^{-sw}\, \frac{{\rm d}w}{2\pi {\rm i}} \nonumber \\
&= \Gamma \bigl( s + \tfrac12 \bigr)^2\, \Pi_h(s)\, J_h(s)\,.
\end{align}

We estimate each of the three factors on the RHS of~\eqref{eq:G_h}
separately. The $\Gamma$-factor and the integral $J_h(s)$
will be estimated everywhere on $\bar D_+(R_h)$, while the product $\Pi_h$
only on the boundary $\partial D_+(R_h)$. By the maximum principle, this will suffice.

We set \[ k(\phi) = \cos\phi + \phi\sin\phi, \quad |\phi|\le\frac{\pi}2\,. \]

\medskip\noindent\underline{Estimate of $J_h(s)$}: We will show that
\[
|J_h(s)| \le \rho^{- 2\rho\cos\phi} e^{2k(\phi)\rho + o_\eps(1) R}, \qquad
s = \rho e^{{\rm i}\phi} \in\bar D_+(R_h)\,.
\]
We start with a standard deformation of the integration contour:
\[
J_h(s) =
\Bigl[\,
\int_{-{\rm i}\pi + 2\log R_h}^{-{\rm i}\pi + 2\log \rho } +
\int_{-{\rm i}\pi + 2\log \rho}^{{\rm i}\pi + 2\log \rho } +
\int_{{\rm i}\pi + 2\log \rho}^{{\rm i}\pi + 2\log R_h }
\,\Bigr]
F_h^\sharp (e^w) e^{\kappa w}\, e^{-sw}\, \frac{{\rm d}w}{2\pi {\rm i}}\,,
\]
and estimate these integrals separately.

First, suppose that $w=u\pm {\rm i}\pi$, $2\log \rho \le u \le 2\log R_h$. Then,
by Lemma~\ref{lemma-SCT},
\[
|F_h^\sharp (e^w) e^{\kappa w}| \le e^{ \frac{h}2 u  + o_\eps(1) R} \cdot e^{-([\frac{h}2] + \frac12)u}
\le e^{o_\eps(1) R}\,,
\]
while $|e^{-sw}| \le  e^{-u\rho\cos\phi + \pi \rho|\sin\phi|}$.
Hence,
\begin{align*}
\Bigl|\,
\int_{\pm {\rm i}\pi + 2\log \rho}^{\pm {\rm i}\pi + 2\log R_h}
F_h^\sharp (e^w) e^{\kappa w}\, e^{-sw}\, \frac{{\rm d}w}{2\pi {\rm i}}
\,\Bigr|
&\le e^{\pi\rho|\sin\phi| + o_\eps(1) R}\,
\int_{2\log\rho}^{2\log R_h} e^{-\rho u \cos\phi\, }\, \frac{{\rm d}u}{2\pi} \\
&< e^{\pi\rho|\sin\phi| + o_\eps(1) R}\,
\int_{\log\rho}^\infty e^{-2\rho u\cos\phi\, }\, \frac{{\rm d}u}{4\pi} \\
&< \rho^{-2\rho \cos\phi}\, e^{\pi\rho|\sin\phi| + o_\eps(1) R}\,.
\end{align*}

Now, consider the case when $w = 2\log\rho+{\rm i}t$, $\pi \le t \le \pi$. In this case,
by Lemma~\ref{lemma-SCT},
$ |F_h^\sharp (e^w) e^{\kappa w} | \le e^{2\rho\cos\frac{t}2 + o_\eps(1) R}$
(since $\kappa=\bigl[\, \tfrac{h}2 \,\bigr]+\tfrac12\, $) and
$ |e^{-sw}| = e^{-2\rho\log\rho\cos\phi +  \rho t\sin\phi} $.
Therefore,
\[
\Bigl|\, \int_{-{\rm i}\pi + 2\log \rho}^{{\rm i}\pi + 2\log \rho }
F_h^\sharp (e^w)e^{\kappa w}\, e^{-sw}\, \frac{{\rm d}w}{2\pi {\rm i}} \,\Bigr|
\le \rho^{-2\rho\cos\phi} e^{o_\eps(1) R}\,
\int_{-\pi}^\pi e^{(2\cos\frac{t}2+t\sin\phi)\rho}\, \frac{{\rm d}t}{2\pi}\,.
\]
Noting that
\begin{equation}\label{eq:k}
\max_{t\in [-\pi, \pi]} \bigl( 2\cos\frac{t}2 + t\sin\phi \bigr) = 2k(\phi), \qquad
|\phi|\le \frac{\pi}2\,,
\end{equation}
we see that
\[
\Bigl|\, \int_{-{\rm i}\pi + 2\log \rho}^{{\rm i}\pi + 2\log \rho }
F_h^\sharp (e^w)e^{\kappa w}\, e^{-sw}\, \frac{{\rm d}w}{2\pi {\rm i}} \,\Bigr|
\le \rho^{-2\rho\cos\phi} e^{2k(\phi)\rho + o_\eps(1) R}\,,
\]
and therefore,
\begin{equation}\label{eq:star}
|J_h(s)| \le \rho^{-2\rho\cos\phi} e^{o_\eps(1)  R}
\Bigl[\, e^{\pi\rho|\sin\phi|} + e^{2k(\phi)\rho} \,\Bigr]\,.
\end{equation}
Letting $t=\pm \pi$ in~\eqref{eq:k}, we get
\[
\pi |\sin\phi| \le 2k(\phi), \qquad |\phi|\le \frac{\pi}2\,,
\]
concluding that the first exponent in the brackets on the RHS of \eqref{eq:star} can be discarded.

\medskip\noindent\underline{Estimate of} $\Pi_h(s)$ on $\partial D_+(R_h)$:
Clearly, $|\Pi_h({\rm i}t)| = 1$. For $s=R_he^{{\rm i}\phi}$, we have
\begin{align*}
\bigl|\, \Pi_h(s) \, \bigr| &= \prod_{1\le \ell \le [\frac{h}2]} \Bigl|\, \frac{1+(\ell - \frac12)/s}{1-(\ell-\frac12)/s}\, \Bigr| \\
&\le \prod_{1\le \ell \le [\frac{h}2]} \, \frac{1+(\ell - \frac12)/R_h}{1-(\ell-\frac12)/R_h}\, \Bigr| \le e^{CR_h^{-1}h^2}\,.
\end{align*}
Hence, everywhere on the boundary of the semi-disk
$D_+(R_h)$, we have \[ \bigl|\, \Pi_h(s) \, \bigr| < e^{\max\{CR_h^{-1}h^2, 1\}}\,. \]

\medskip\noindent\underline{Estimate of} $\Gamma \bigl( s+\tfrac12 \bigr)^2$:
By the Stirling formula,
\[
\Gamma \bigl( s+\tfrac12 \bigr) =
\sqrt{2\pi}\, s^s e^{-s} (1+o(1)), \quad \text{Re\,}s\to +\infty\,,
\]
whence,
\begin{align*}
\bigl| \Gamma \bigl( s+\tfrac12 \bigr) \bigr|^2 \lesssim
\rho^{2\rho\cos\phi}\, e^{-2\rho(\phi\sin\phi +\cos\phi)}
= \rho^{2\rho\cos\phi}\, e^{-2k(\phi)\rho}\,, \qquad {\rm Re}(s)\ge 0\,.
\end{align*}

\medskip Combining the estimates of $J_h(s)$, $\Pi_h(s)$, and $\Gamma\bigl(s+\tfrac12 \bigr)^2$,
we see that $|g_h(\la)| < e^{o_\eps(1) R +O(R^{-1}h^2)}$ everywhere
on the boundary $\partial \Omega_h$, and hence, inside $\Omega_h$, as well.  \hfill $\Box$

\subsection{Some estimates}

For an interval $I\subset\mathbb R$, we denote by $L(I)$ its length, and by $I_{+t}$ its open complex
$t$-neigbourhood.

The following Jensen-type bound is classical.
To keep the exposition self-contained, we recall its proof.

\begin{lemma}\label{lemma:bounds0}
Let $I\subset \bR$ be an interval of length
\begin{equation}\label{eq:L}
L(I) \lesssim t,
\end{equation}
and let $\Lambda\subset I$ be a discrete set of cardinality $|\Lambda|$.
Then there exists a positive constant $c$,
depending only on the implicit constant in~\eqref{eq:L}, such
that, for any function $f$ analytic in $I_{+t}$ and vanishing on $\Lambda$, we have
\[
\max_{\overline{I}_{+t/2}} |f| < e^{-c|\Lambda|} \sup_{I_{+t}} |f|\,.
\]
\end{lemma}

\noindent{\em Proof}:  Consider the Riemann mapping
$\phi\colon \bD\to I_{+t}$ which maps the origin to the center of the interval $I$.
We claim that there exists $\rho<1$, independent of $t$, such that
$\overline{I}_{+t/2} \subset \phi (\rho\bD)$. It suffices to check this only for $t=1$
(otherwise, we replace the function $\phi$ by $t\phi$), so we assume that $t=1$. If the claim
does not hold, then there exists a sequence of points $(z_n)\subset \bD$ such that $z_n\to\partial \bD$, but $\phi(z_n)\in \overline{I}_{+1/2}$, which contradicts the covering property of Riemann mappings.

Now, consider the function $ f \circ \phi$. It vanishes at the points
of $\phi^{-1}(\Lambda)\subset\rho\bD$. Denote by $B_a(z)=(z-a)/(1-z\bar a)$ the Blaschke factor,
and set \[ c(\rho) = \max \{ |B_a(z)|\colon |a|, |z| \le \rho\}. \]
Obviously, $c(\rho)<1$. Then
$ f\circ \phi = g \cdot \prod_{a\in\phi^{-1}(\Lambda)} B_a $, where the function $g$
is analytic in $\bD$, and $\sup_\bD |g| = \sup_{\bD} |f\circ\phi|$. Hence,
\begin{align*}
\max_{\bar I_{+z/2}} |f| \le \max_{\rho\bar\bD} |f\circ \phi|
&\le c(\rho)^{|\Lambda|} \max_{\rho\bar\bD} |g| \\
&\le e^{-c|\Lambda|} \sup_{\bD} |g| = e^{-c|\Lambda|} \sup_{\bD} |f\circ \phi | =
e^{-c|\Lambda|} \sup_{I_{+t}} |f|,
\end{align*}
proving the lemma. \hfill $\Box$

\begin{lemma}\label{lemma:bounds1}
Let $t\ge t_0$. Let $I\subset \bR$ be an interval of length
\begin{equation}\label{eq:L1}
t_0 \le L(I) \lesssim t,
\end{equation}
and let $\Lambda\subset I$ be a discrete set of cardinality
\begin{equation}\label{eq:La1}
|\Lambda|\gtrsim t.
\end{equation}
Then there exists
a positive constant  $c$ (depending only on the implicit constants in~\eqref{eq:L1} and~\eqref{eq:La1}) such that, for any function $g$ analytic function in $I_{+t}$ satisfying
\begin{equation}\label{eq0}
\sup_{I_{+t}} |g| < e^{o_\eps(1) t},
\end{equation}
and
\[
|g(\la)|=1, \qquad \la\in\Lambda,
\]
we have
\begin{align}
\label{eq1} & \max_{s\in I} ||g(s)|-1| < e^{-ct}, \\
\label{eq2} & \max_{s\in I}\, \Bigl|\, \frac{g(s+1)}{g(s)} - 1\Bigr| < o_\eps (1) + e^{-ct},
\intertext{and}
\label{eq3} 
& \max_{s\in I}\, \Bigl|\, \frac{g(s)^2}{g(s-1)g(s+1)} - 1\Bigr|
< o_\eps (1) + e^{-ct},
\end{align}
provided that $\eps \ll 1 \ll t_0$.
\end{lemma}

\noindent{\em Proof}: Set $g^*(s)=\overline{g(\bar s)}$ and $f(s)= g(s)g^*(s)-1$.
The function $f$ is analytic on $I_{+t}$, satisfies~\eqref{eq0}, and vanishes on $\Lambda$.
Hence, by Lemma~\ref{lemma:bounds0},
\[
\max_{\overline{I}_{+t/2}} |f| < e^{-ct} + e^{o_\eps(1)t} < e^{-ct}\,.
\]
In particular,
\[
\max_{I} |\,|g|-1| \le \max_I |\, |g|^2-1| = \max_I |f| < e^{-ct}\,,
\]
proving~\eqref{eq1}.

Set $I^* = I_{+1}\cap \bR$. For $s\in I$, we have
\[
\frac{g(s+1)}{g(s)} = \exp\Bigl( \int_s^{s+1} \frac{g'}g(\xi)\, {\rm d}\xi \Bigr),
\]
whence,
\[
\Bigl| \frac{g(s+1)}{g(s)} - 1 \Bigr| \le 2\max_{I^*} |g'/g|\,,
\]
assuming that the maximum on the RHS is sufficiently small.

Take $\xi\in I^*$, and let $\rho = \tfrac12 t$.
Since the function $f$ is very small on $\overline{I}_{+t/2}$, the function $g$
cannot vanish therein, and therefore, we can choose in $I_{+t/2}$ a branch of $\log g$.
Applying the Schwarz formula to $\log g$ in the disk $D(\xi, \rho)\subset I_{+t/2}$,
we have
\begin{equation}
\label{eq:Poisson}
\log g(\xi+s) = \frac1{2\pi}\,
\int_{-\pi}^\pi \log |g(\xi+\rho e^{{\rm i}\theta})|\,
\frac{\rho e^{{\rm i}\theta}+s}{\rho e^{{\rm i}\theta}-s}\, {\rm d}\theta + {\rm i}C.
\end{equation}
Differentiating in $s$ and letting $s=0$, we get
\[
\frac{g'}g(\xi) = \frac1{\pi\rho}\,
\int_{-\pi}^\pi \log |g(\xi+\rho e^{{\rm i}\theta})|\, e^{-{\rm i}\theta}\, {\rm d}\theta\,,
\]
whence
\begin{equation}\label{eq:g'/g}
\Bigl|\, \frac{g'}g(\xi)\, \Bigr| \le \frac1{\pi\rho}\,
\int_{-\pi}^\pi \bigl| \log |g(\xi+\rho e^{{\rm i}\theta})|\bigr|\, {\rm d}\theta\,.
\end{equation}

Everywhere in the disk $D(\xi, \rho)$, we have $\log |g| \le o_\eps (1) t$, while in the center of
that disk $\log |g|$ is very small:
\[
|\log |g(\xi)|\,| = |\log(1 + (|g(\xi)|-1))| < e^{-ct}.
\]
Recalling that $\log |g|$ is harmonic in the disk $D(\xi, \rho)$, we get
\begin{multline*}
\Bigl| \int_{-\pi}^\pi \bigl[ \log_+ |g(\xi+\rho e^{{\rm i}\theta})|
- \log_- |g(\xi+\rho e^{{\rm i}\theta})| \bigr] \, {\rm d}\theta \Bigr| \\
= \Bigl| \int_{-\pi}^\pi \log |g(\xi+\rho e^{{\rm i}\theta})|\, {\rm d}\theta \Bigr|
= |\log |g(\xi)|\,| < e^{-ct}.
\end{multline*}
Therefore,
\[
\int_{-\pi}^\pi \log_- |g(\xi+\rho e^{{\rm i}\theta})|\, {\rm d}\theta
< \int_{-\pi}^\pi \log_+ |g(\xi+\rho e^{{\rm i}\theta})|\, {\rm d}\theta + e^{-ct}
< o_\eps(1)t + e^{-ct}\,,
\]
and
\begin{multline*}
\int_{-\pi}^\pi |\log |g(\xi+\rho e^{{\rm i}\theta})|\,|\, {\rm d}\theta \\
= \int_{-\pi}^\pi \bigl[ \log_+ |g(\xi+\rho e^{{\rm i}\theta})|\, {\rm d}\theta
+ \log_- |g(\xi+\rho e^{{\rm i}\theta})| \bigr]\, {\rm d}\theta <  o_\eps(1)t + e^{-ct}\,.
\end{multline*}
Plugging this estimate
in~\eqref{eq:g'/g} and recalling that $\rho$ and $t$ are comparable,
we conclude that $\max_{I^*} |g'/g| <  o_\eps(1) + e^{-ct}$, proving~\eqref{eq2}.

At last,
\[
\Bigl|\, \frac{g(s)^2}{g(s-1)g(s+1)} - 1\Bigr| \le 2\max_{I^*} |(g'/g)'|,
\]
providing that the maximum on the RHS is sufficiently small.
Differentiating twice the Schwarz formula~\eqref{eq:Poisson}
and letting $s=0$, we see that
\[
\Bigl( \frac{g'}{g}\Bigr)'(\xi) = \frac2{\pi\rho^2}\,
\int_{-\pi}^\pi \log|g(\xi+\rho e^{{\rm i}\theta})|
e^{-2{\rm i}\theta}\, {\rm d}\theta,
\]
whence
\[
\max_{I^*} |(g'/g)'|
\le \frac2{\pi\rho^2}\, \int_{-\pi}^\pi \bigl| \log|g(\xi+\rho e^{{\rm i}\theta})| \bigr|\,
{\rm d}\theta
\lesssim t^{-1}o_\eps(1) + e^{-ct},
\]
proving~\eqref{eq3}. \hfill $\Box$

\subsection{Local bounds of $F$}
We assume that $F$ is an {\sf EFET}$(\sigma)$ satisfying assumptions of Theorem~\ref{thm:main}, that is,
\[
F(z) = \sum_{n\ge 0} \omega_n\, \frac{z^n}{n!}\,,
\]
where the set $\Lambda_F =\{n\in\bZ_+\colon |\omega_n|=1$\} has a lower density $d>0$.
Here we will show that if $F$ has relatively few zeroes, then it has bounds~\eqref{eq_G} and~\eqref{eq_alpha} needed for application of Lemma~\ref{lemma_interpol}.

\begin{lemma}\label{lemma:local_bounds}
Let $F$ satisfy the assumptions of Theorem~\ref{thm:main}, and let
there exist 

\smallskip\noindent$\bullet$ a sufficiently small $\eps$, depending only on the type $\sigma$ of $F$
and on the lower density~$d$ of the set $\Lambda_F$, 

\smallskip\noindent$\bullet$  $A=\omega_\eps(1)$,

\smallskip\noindent$\bullet$ 
and a sufficiently large $R=R_\eps$, 

\smallskip\noindent such that
\[
n_F(AR) < \eps R.
\]
Then there exist 

\smallskip\noindent$\bullet$ 
$\alpha\in \mathbb C$, $|\alpha| < 1 + o_\eps(1)$, 

\smallskip\noindent$\bullet$ and $A_1 = \omega_\eps(1)$

\smallskip\noindent such that 
\[
F(z) = e^{\alpha z} G(z),
\]
with \[\log M_G(A_1R) < o_\eps (1) R.\]
\end{lemma}
\noindent{\em Proof of Lemma~\ref{lemma:local_bounds}} will consist of several steps.
We fix $A_1  = \omega_\eps(1)$ so that \[ A_1^2 \max\{\eps, A^{-1}\} = o_\eps(1)\,. \]

\medskip\noindent\underline{Bounding $G$}:
We choose $r_0$ so large that
$\log M_F(r) \le (\sigma + \eps)r$ and $n_F(r)\lesssim \sigma r$ for $r\ge r_0$,
and assume that $r_0$ is a way smaller than $R$. We let $P$ be a polynomial
vanishing at zeroes of $F$ lying in the disk $\{|\la|\le r_0\}$ and normalized
so that $P(0)=F(0)$.

Set $R_1=A_1 R$. Applying the Hadamard factorization theorem to $F(z)$ with $|z|=R_1$, we get
\begin{align*}
F(z) &= e^{\alpha z} P(z) \prod_{r_0 < |\la| \le AR} \Bigl( 1 - \frac{z}\la \Bigr)
\prod_{|\la| > AR } \Bigl( 1 - \frac{z}\la \Bigr)e^{\frac{z}{\la} } \\
&= e^{\alpha z} P(z) \pi_1 (z) \pi_2(z) = e^{\alpha z} G(z).
\end{align*}

We have
$ \log |P(z)| \le C_F + n_F(r_0) \log (R+r_0) = o_\eps(1) R $,
provided that $R$ is sufficiently large. We also have
\[
\log |\pi_1 (z)| \le \int_{r_0}^{AR} \log\Bigl( 1 + \frac{R_1}t \Bigr)\, {\rm d}n_F(t)\,.
\]
Integrating by parts, we bound the integral on the RHS by
\begin{align*}
n_F(AR)&\log\Bigl( 1 + \frac{A_1}{A} \Bigr) + R_1\, \int_{r_0}^{AR} \frac{n_F(t)}{t(t+R_1)}\, {\rm d}t
\\
&\le R_1\, \int_{r_0}^{AR} \frac{n_F(t)}{t(t+R_1)}\, {\rm d}t + o_\eps(1) R \\
&< \int_{r_0}^{R/A_1} \frac{n_F(t)}t\, {\rm d}t
+ R_1 \int_{R/A_1}^{AR} \frac{n_F(t)}{t^2}\, {\rm d}t + o_\eps(1) R \\
&\lesssim \int_{r_0}^{R/A_1} \sigma \, {\rm d}t
+ R_1 \int_{R/A_1}^{AR} \frac{\eps R}{t^2}\, {\rm d}t + o_\eps(1) R
\\
&\lesssim \bigl( \sigma A_1^{-1} + A_1^2 \eps + o_\eps(1) \bigr) R = o_\eps(1) R,
\end{align*}
and similarly
\begin{align*}
\log |\pi_2 (z)| &\lesssim R_1^2\, \int_{AR}^\infty \frac{{\rm d}n_F(t)}{t^2} \\
&< 2R_1^2\, \int_{AR}^\infty \frac{n_F(t)}{t^3}\, {\rm d}t \\
&\lesssim (A_1R)^2\, \int_{AR}^\infty \sigma\, \frac{{\rm d}t}{t^2} \\
&\lesssim  \sigma A_1^2 A^{-1} R = o_\eps(1) R\,.
\end{align*}
Combining these bounds, we see that $\log M_G(A_1R) = o_\eps(1)R$.

\medskip\noindent\underline{Preliminary bound for $|\alpha|$}: Here we show that
\[
|\alpha| < R^{-1}\, \log M_F(R) + o_\eps(1)\,.
\]
We assume that $\alpha$ is positive (otherwise, we rotate $z$) and
that $F(0)\ne 0$, the case $F(0)=0$ does not bring new difficulties. Then we  take a sufficiently small $\kappa>0$ (to be chosen momentarily),
and write
\begin{align*}
 -C_F = \log |F(0)| &= \log |G(0)| \\
&\le \int_{-\pi}^\pi \log |G(Re^{{\rm i}\theta})|\, \frac{{\rm d}\theta}{2\pi} \\
&= \int_{|\theta|\le \kappa} \bigl[ \log |F(Re^{{\rm i}\theta})|
- \alpha  R \cos\theta \bigr]\, \frac{{\rm d}\theta}{2\pi} +
\int_{\kappa \le |\theta| \le \pi } \log |G(Re^{{\rm i}\theta})|\, \frac{{\rm d}\theta}{2\pi} \\
&\le 2\kappa \log M_F(R) - \alpha R (2\kappa - \kappa^3/3)) + o_\eps(1) R\,,
\end{align*}
whence,
\begin{align*}
\alpha &\le
(1+C\kappa^2)\, R^{-1}\log M_F(R) + \kappa^{-1} o_\eps(1) + C_F(\kappa R)^{-1}
\\
&\le
R^{-1}\log M_F(R) + C \bigl( \sigma\kappa^2 + \kappa^{-1}o_\eps(1) \bigr) + o_\eps(1).
\end{align*}
Choosing $\kappa = o_\eps (1)$ so that $(\sigma \kappa^2 + \kappa^{-1} o_\eps(1)) = o_\eps (1)$, 
we see that the RHS is bounded by \[ R^{-1}\log M_F(R) + o_\eps (1), \]
provided that $R$ is sufficiently big.

\medskip\noindent\underline{Bounding $|\omega_n|$}: Let $I=[bR, BR]$ where parameters
$B$ and $b$, $0<b\ll 1 \ll B$ (depending only on the type $\sigma$ and density $d$) are to be chosen.
Here we show that, for $n\in I$, the absolute values of the coefficients
$|\omega_n|$ cannot be essentially bigger than $1$. To show this, we will use 
the analytic interpolation Lemma~\ref{lemma_interpol}.

We set
\[
\widetilde{F}(z) = F\Bigl( \frac{z}{\alpha}\Bigr)
= \sum_{n\ge 0} \frac{\omega_n}{\alpha^n}\, \frac{z^n}{n!}\,.
\]
Then $n_{\widetilde F} (AR/\alpha) < o_\eps(1) R$, and, by the previous steps,
\[ \widetilde{F}(z) = e^z G\Bigl( \frac{z}{\alpha}\Bigr) = e^z \widetilde{G}(z) \]
with
\[
\log M_{\widetilde{G}} (r) = \log M_G (r/\alpha) < o_\eps (1)R,
\]
whenever $r<\alpha A_1 R$, in particular, for $r=2BR+\tfrac12$.
This allows us to apply Lemma~\ref{lemma_interpol} to the function $\widetilde F$ with $h=0$ (and $2BR+\tfrac12$ instead of $R$).
It gives us the interpolation function $g_0$,
\[
g_0(n)=\Bigl|\, \frac{\omega_n}{\alpha^n}\, \Bigr|^2, \qquad n\in\bZ_+\,,
\]
analytic in the right half-plane and satisfying
\[
|g_0(s)| < e^{o_\eps(1) R}
\]
in the closed semi-disk $\bar \Omega = \{s\colon |s|\le 2BR, \text{\ Re}(s) \ge 0\}$.

Next, we set
\[
f(s) = g_0(s) - \frac1{|\alpha|^{2s}}\,.
\]
Then the function $f$ has the same upper bound in $\bar\Omega$ as $g_0$, i.e.,
\[
|f(s)| < e^{o_\eps(1) R}, \quad z\in\bar\Omega\,,
\]
and
\[
f(n) = \Bigl| \frac{\omega_n}{\alpha^n} \Bigr|^2 - \frac1{|\alpha|^{2n}}\,.
\]
In particular, $f(n)=0$ for $n\in\Lambda_F$.

We set $t=bR$. Then $I_{+t}\subset\Omega$.
We assume that $b<d/2$ and $B\ge 2$. Then
\[
|\Lambda_F \cap I| \ge \frac{d}3\, BR = c_0,
\]
and, by Lemma~\ref{lemma:bounds0},
\[
\max_{\bar I_{+t/2}} |f| \le e^{-cR} \sup_{I_{+t}} |f|\,.
\]
Thus,
\[
|f(s)| < e^{-cR +o_\eps(1)R} < e^{-c_1 R}, \qquad s\in \overline{I}_{+t/2},
\]
provided that $\eps$ is sufficiently small. In particular this estimate holds on $I$, and therefore,
\[
\Bigl| \frac{\omega_n}{\alpha^n} \Bigr|^2 - \frac1{|\alpha|^{2n}} = f(n) <
e^{-c_1 R}\,.
\]
Hence,
\begin{align*}
|\omega_n| &< \bigl( 1 + |\alpha|^{2n} e^{-c_1 R} \bigr)^{1/2} \\
&\le \sqrt{2}\, \max\bigl\{ 1, (|\alpha|e^{-c_2})^n\, \bigr\}, \qquad n\in I,
\end{align*}
completing this step.

\medskip\noindent\underline{Bounding $F$}:
Now, we are ready to bound $|F(z)|$ for $|z|=R$.
We have
\[
|F(z)| \le \Bigl( \sum_{0\le n<B^{-1}R} + \sum_{B^{-1}R \le n \le BR} + \sum_{n>BR}\Bigr) \frac{|\omega_n|}{n!}\, R^n.
\]
Since $F(z)$ has finite exponential type, we have $|\omega_n| \le C Q^n$ with $Q$ depending on the type of $F$ (for instance, we can take $Q=2\sigma$). Hence, assuming that $b<B^{-1}$, i.e., that
$[B^{-1}R, BR]\subset I$, we have
\[
|F(z)| \lesssim
R \frac{(QR)^{B^{-1}R}}{[B^{-1}R]!}
+ \sum_{B^{-1}R \le n \le BR} \frac{\max\bigl\{ 1, ( |\alpha|e^{-c_2})^n \bigr\}}{n!}\, R^n
+ \frac{(QR)^{BR}}{([BR]+1)!}.
\]
The first term on the RHS is bounded by
\[
R \bigl( QBe \bigr)^{B^{-1}R} < e^{R},
\]
provided that $B^{-1}$ is sufficiently small. Similarly, the third term on the RHS  does not exceed
\[
\Bigl( \frac{Qe}{B} \Bigr)^{BR} < 1,
\]
provided that $B$ is sufficiently large. At last,
\[
\sum_{B^{-1}R \le n \le BR} \frac{\max\bigl\{ 1, ( |\alpha|e^{-c_2})^n \bigr\}}{n!}\, R^n
< \max\bigl\{ e^{R}, e^{|\alpha|e^{-c_2}R}\bigr\},
\]
whence,
\[
\log M_F(R) < \max\{ 1,  |\alpha|e^{-c_2} \} R + O(1),
\]

\medskip\noindent\underline{Completing the proof of Lemma~\ref{lemma:local_bounds}}:
\begin{align*}
|\alpha| &< R^{-1}\, \log M_F(R) + o_\eps(1) \\
&< \max\{1, |\alpha|e^{-c_2}\} + o_\eps(1)\,,
\end{align*}
which yields that $|\alpha|\le 1+o_\eps(1)$. \hfill $\Box$

\subsection{A combinatorial lemma}

\begin{lemma}\label{lemma:combi}
Let $\Lambda\subset \bZ_+$ be a set of positive lower density $d$.
Then there exist positive $c_1<\tfrac12$ and $c_2$ such that
given sufficiently big $R>1$,
one can find integers $x$, $2c_1 R < x < (1-2c_1) R$, and $h$, $c_1R \le h \le (1-c_1)R$,
satisfying the following:

Let $J=\Lambda\cap [1, c_1R]$ and $K=\Lambda\cap [x, x+c_1 R]$. Then
\[ |K|\ge \frac{d}2\, c_1R,  \quad |(J+h)\cap K| \ge c_2 R.\]
\end{lemma}
\noindent Note that in the assumptions of this lemma we also have
\[
|J| \ge \frac{d}2\, c_1 R
\]
(since we assume that $R$ is sufficiently big).

\medskip\noindent{\em Proof of Lemma~\ref{lemma:combi}}: We split $[0, R]$ into segments of length
$c_1 R$. The two first and two last segments contain together at most
$\displaystyle 4c_1 R < \frac{d}2\, R$ elements of $\Lambda$ (we choose $c_1$ sufficiently small), so the rest of the segments contains at least $\displaystyle \frac{d}2\,R$ elements of $\Lambda$. Hence, at least one of the segments contain not less than $\displaystyle \frac{c_1d}2\,R$ elements of $\Lambda$.
We choose this segment, denote it by $[x, x+c_1 R]$ and use it to define the set $K$.

To find $h$, let us look at a translation of $J$ by a uniformly chosen $h'$ in the range
$x-c_1 R<h'<x+c_1 R$. For each $j\in J$, the chance that $j+h'\in K$ is at least
\[ \frac1{2c_1R}\, |K| \ge \frac{d}5\,, \] so
\[
\bE\bigl[ |(J+h')\cap K|\bigr] \ge \frac{d}5\, |J| \ge \frac{d^2}{10}\, c_1 R = c_2 R\,.
\]
In particular, there exists a specific $h$ for which $|(J+h)\cap K|$ achieves the expected value.
\hfill $\Box$

\subsection{All coefficients $\omega_n$ have absolute value $1$}
We will prove Theorem~\ref{thm:main} in two steps.
Here, assuming that $F$ is an entire function of
exponential type such that, for a sequence $R_j\uparrow \infty$,
\begin{equation}\label{eq:lower_limit}
n_F(R_j)<\eps R_j,
\end{equation}
with sufficiently small $\eps$ depending only on the exponential type of $F$ and on the lower density of the set $\Lambda_F$, we show that, for all $n\in\bZ_+$, $|\omega_n|=1$.
Changing $\eps$, we turn~\eqref{eq:lower_limit} into a seemingly stronger condition that
\[ n_F(AR_j)<\eps R_j\,, \] with some $A=\omega_\eps(1)$. For instance, we can take $A=1/\sqrt \eps$,
and then replace $R_j$ by $R_j/A$ and $\eps$ by $\sqrt\eps$.

We fix a sufficiently large $R=R_j$ depending on $\eps$.
Applying Lemma~\ref{lemma:local_bounds}, we factorize $F(z)=e^{\alpha z}G(z)$, with $\log M_G(A_1R)<o_\eps(1)R$ and $|\alpha|<1+o_\eps(1)$. Then, given $h\in\bZ_+$,
Lemma~\ref{lemma_interpol}, applied with $A_1R$ (instead of $R$),
provides us with an interpolating function $g_h$, $g_h(n)=\omega_n \bar\omega_{n+h}$, holomorphic in the closed semi-disk
\[
\bar \Omega_h=\{\la\colon \text{Re}(\la)\ge -\kappa, |\la+\kappa|\ge A_1R-\kappa\}, \quad \kappa=[h/2]+1/2,
\]
and satisfying
\[
\max_{\bar\Omega_h} |g_h| \le e^{o_\eps(1)A_1 R +O(h^2/(A_1R))}\,.
\]
Replacing, if needed, $A_1$ by a smaller quantity which is still $\omega_\eps(1)$, we will assume
that
\[
\max_{\bar\Omega_h} |g_h| \le e^{o_\eps (1)(R+R^{-1}h^2)}\,.
\]

First we look at the interpolating function $g_0$, which satisfies
\[
|g_0(\la)| < e^{o_\eps(1)R}, \qquad \la\in\bar\Omega_0,
\]
and
\[
g_0(n) = |\omega_n|^2, \qquad n\in \bZ_+.
\]
In particular, $g_0(n)=1$ for $n\in \Lambda_F$.
We apply Lemma~\ref{lemma:bounds1} to the function $g_0$ and the interval $I=[cR, R]$
with $t=cR$ (note that $[cR, R]_{+t}\subset\Omega_0$).
Choosing sufficiently small $c<d/2$ and letting $R$ be sufficiently large, we ensure that
\[
|\Lambda_F \cap [cR, R]| \ge cR,
\]
which gives us good bounds in the range $cR\le n \le R$:
\[
|\omega_n| = 1 + e^{-c_1 R}
\]
and similarly
\[
\Bigl| \frac{\omega_{n+1}^2}{\omega_n \omega_{n+2}} \Bigr| = 1 + O(R^{-1})\,.
\]

To descend to small values of $n$, we apply the combinatorial lemma,
which provides us with
the sets $J$ and $K$ and the integer $h$ satisfying:
\begin{itemize}
\item $J=\Lambda_F\cap [0, c_1R]$, $|J|\ge c_2R$;
\item $K=\Lambda_F\cap [x, x+c_1R]$, where $2c_1R<x<(1-2c_1)R$, and $|K|\ge c_2R$;
\item $|J\cap (K-h)|\ge c_2R$;
\item $c_1R \le h \le (1-c_1)R$
\end{itemize}
(if necessary, we increase $R$), and consider the interpolating function $g_h$ with this value of the parameter $h$.
The absolute value of $g_h(n)=\omega_n\bar\omega_{n+h}$ equals $1$ on a sufficiently fat set
$J \cap (K-h) \subset [0, c_1R]$, and the $\tfrac12\, c_1R$-neighbourhood of the interval
$[0, c_1R]$ lies in the semi-disk $\Omega_h$. Hence, by the second quotient bound from Lemma~\ref{lemma:bounds1}
applied to $g_h$ with $I=[0, c_1R]$ and $t= c_1R/2$, we have
\[
\frac{g_h(n)g_h(n+2)}{g_h(n+1)^2} = 1 + O(R^{-1}),
\]
that is,
\[
\frac{\omega_n \omega_{n+2}}{\omega_{n+1}^2} \cdot
\frac{\bar\omega_{n+h} \bar\omega_{n+h+2}}{\bar\omega_{n+h+1}^2}
= 1 + O(R^{-1}).
\]
But, for $n\in [0, cR]$ and $c_1 R \le h \le (1-c_1)R$, $c<c_1$,
we have \[ n+h, n+h+1, n+h+2\in [c_1R, R], \]
so
\[
\Bigl| \frac{\bar\omega_{n+h} \bar\omega_{n+h+2}}{\bar\omega_{n+h+1}^2}\, \Bigr|
= 1 + O(R^{-1})
\]
and therefore,
\[
\Bigl| \frac{\omega_n \omega_{n+2}}{\omega_{n+1}^2} \Bigr|
= 1 + O(R^{-1}).
\]
Letting $R\to\infty$ along a sequence of values where~\eqref{eq:lower_limit} holds,
we get
\[
\Bigl| \frac{\omega_n \omega_{n+2}}{\omega_{n+1}^2} \Bigr| = 1
\]
for all $n\in\bZ_+$.

Thus, the quotient $|\omega_n|/|\omega_{n+1}|$ is constant for all $n\in\bZ_+$, that is,
$|\omega_n| = \beta \gamma^n$. We know that, for every $n\in\Lambda_F$, $|\omega_n|=1$,
so $\beta=\gamma=1$. This completes the first step of the proof. \hfill $\Box$

\subsection{Coup de gr\^ace}

Here, our strategy will be similar to the one used in the first step.
Since we already know that all the coefficients $\omega_n$ are unimodular, we will freely
use that
\[
g_h(n) = \omega_n \bar\omega_{n+h} = \frac{\omega_n}{\omega_{n+h}}\,.
\]

First, we look at the interpolating function $g_1$. It takes unimodular values at the non-negative integers
and its absolute value is bounded by $e^{o_\eps(1)R}$ on a $cR$-neighbourhood of the interval $[cR, R]$. Hence, by Lemma~\ref{lemma:bounds1} applied with $t=cR$,
\begin{equation}\label{eq:double_quotient}
\frac{\omega_{n+1}^3 \omega_{n+3}}{\omega_n\omega_{n+2}^3}
= \frac{g_1(n+1)^2}{g_1(n)g_1(n+2)} = 1 + O(R^{-1})\,, \qquad n\in [cR, R]\,.
\end{equation}

Then, we look at the function $g_h$ with $h=[R]$. It also takes unimodular values at non-negative integers, and
its absolute value  is
bounded by $e^{o_\eps(1)R}$ on a $cR$-neighbourhood of the interval $[0, cR]$. Applying again Lemma~\ref{lemma:bounds1} with $t=cR$,
we deduce that
\begin{multline*}
\frac{\omega_{n+1}^3 \omega_{n+3}}{\omega_n \omega_{n+2}^3}\, :\,
\frac{\omega_{n+h+1}^3 \omega_{n+h+3}}{\omega_{n+h}\omega_{n+h+2}^3}
= \frac{g_h(n+1)^2}{g_h(n)g_h(n+2)}\, :\, \frac{g_h(n+2)^2}{g_h(n+1)g_h(n+3)} \\
= 1 + O(R^{-1})\,,  \qquad n\in [0, cR]\,.
\end{multline*}
But now $n+h\in [cR, R]$, so that, by~\eqref{eq:double_quotient},
\[
\frac{\omega_{n+h+1}^3 \omega_{n+h+3}}{\omega_{n+h}\omega_{n+h+2}^3} = 1 + O(R^{-1}),
\]
and therefore,
\[
\frac{\omega_{n+1}^3 \omega_{n+3}}{\omega_n \omega_{n+2}^3} = 1 + O(R^{-1})\,, \qquad
n\in [0, cR]\,.
\]
Letting $R = R_j\uparrow\infty$,
we conclude that
\[
\frac{\omega_{n+1}^2}{\omega_n\omega_{n+2}}\, :\, \frac{\omega_{n+2}^2}{\omega_{n+1}\omega_{n+3}}
= \frac{\omega_{n+1}^3 \omega_{n+3}}{\omega_n \omega_{n+2}^3} = 1, \qquad n\in\bZ_+\,,
\]
which means that the double quotient
\[
\frac{\omega_{n+2}}{\omega_{n+1}}\,:\, \frac{\omega_{n+1}}{\omega_{n}}
\]
is constant, so
\[
\omega_{n+1} = C^{n(n+1)/2}\, \frac{\omega_1}{\omega_0}\,.
\]
Thus, $\omega_n = e^{2\pi {\rm i}(\beta n^2 + \gamma n + \delta)}$ for some
$\beta, \gamma,  \delta \in [0, 1)$.

\medskip
To complete the proof,
we claim that $\beta=0$ or $\beta=\tfrac12$. To show this, we consider
the interpolating function $g_d$ with bounded $d>0$. Applying Lemma~\ref{lemma:bounds1}
on the interval $[cR, R]$, we see that
\begin{equation}\label{eq:g_d}
\frac{g_d(n+1)}{g_d(n)} = 1 + o_\eps(1) + e^{-c_1R}\,,  \qquad n\in [cR, R]\,.
\end{equation}
On the other hand,
\begin{align*}
\frac{g_d(n+1)}{g_d(n)} &= \frac{\omega_{n+1} \omega_{n+d}}{\omega_n \omega_{n+d+1}} \\
&= e^{2\pi {\rm i}\beta ((n+1)^2 + (n+d)^2 - n^2 - (n+d+1)^2)}
= e^{-2\pi {\rm i} (2\beta d)}\,.
\end{align*}
If $2\beta$ is not an integer, we can find $d=O(\beta^{-1})$ such that
\[
\frac14 < 2\beta d < \frac34 \qquad (\text{mod\ }1),
\]
which yields contradiction with~\eqref{eq:g_d}
when $R$ is sufficiently large and $\eps$ is small enough.

At last, if $\beta = \tfrac12$, we have
\[
e^{2\pi {\rm i} \beta n^2} = e^{\pi {\rm i} n^2} = e^{\pi {\rm i} n}\,,
\]
so we can let $\beta = 0$ adding $\tfrac12$ to $\gamma$. This finishes off the proof of the
theorem. \hfill $\Box$

\medskip\noindent{\bf Remark:}
At the very last step, we could complete the proof referring to a result of Eremenko and Ostrovskii~\cite[Theorem~3]{EO}
on the zero distribution of entire functions of the form
\[
\sum_{n\ge 0} e^{2\pi{\rm i}\beta n^2}\, \frac{z^n}{n!}
\]
(for a different proof of a more general result see~\cite[Theorem~1]{BNS}).
Since what we used is a simple consequence of the technique developed above, we decided to follow a self-contained argument.

\section{Proof of Theorem~\ref{thm:main2}}

Here, 
\[
F(z) = \sum_{n\ge 0} \omega_n\, \frac{z^n}{n!}
\]
with $0 < c \le |\omega_n| \le C <\infty$, $n\in\mathbb Z_+$.
WLOG, we assume that $F(0)=1$. 

First, we notice that in the assumptions of the theorem, we have
\begin{equation}\label{eq:10}
R-\frac12 \log R - O(1) \le \log M_F(R) < R + O(1), \qquad R\to\infty.
\end{equation}
The upper bound is obvious. To prove the lower bound it suffices to assume that $R$ is a 
large positive integer.
Then
\[
M_F(R)^2 \ge \int_{-\pi}^\pi |F(Re^{{\rm i}\theta})|^2\, \frac{{\rm d}\theta}{2\pi} 
= \sum_{n\ge 0} \frac{|\omega_n|^2}{(n!)^2}\, R^{2n} \gtrsim \frac{R^{2R}}{(R!)^2}
\gtrsim e^{2R - \log R}\,.
\]

\subsection{Proof of part (i)}
By Hadamard's factorization theorem, we have
$F(z) = e^{az} G(z)$, where $G$ is a canonical product of genus $0$. 
Since $n_G(R) = o(\sqrt{R})$ as $R\to\infty$, $G$ has minimal type 
with respect to the order $\tfrac12$. Then
$\log M_F(r)=|a|r + o(\sqrt{r})$, and by~\eqref{eq:10}, $|a|=1$. To simplify notation, we assume 
that $a=1$. Then
\[
|G(R)| = e^{-R}|F(R)| \le e^{-R}M_F(R) = O(1)\,,
\]
and, by one of the classical versions of the Phragm\'en--Lindel\"of theorem~\cite[Lecture~6, Theorem~2]{Levin}, 
$G$ is a constant function, that is, $F(x) = C e^{z}$. \hfill $\Box$

\medskip\noindent{\bf Remark:} It is not difficult to see that the very similar proof  
yields the same conclusion under less stringent bounds on the sequence of multipliers $\omega_n$.
It suffices to assume that
\[ |\omega_n|=o(n), \quad n\to\infty, \] 
and 
\[
 \limsup_{n\to\infty} |\omega_n|^{1/n} = 1\,.
\]

\subsection{Proof of part (ii)}

Now, we assume that 
\[
\liminf_{R\to\infty} \frac{n_F(R)}{R^\alpha} = 0
\]
with some $\alpha<\tfrac12$. We fix $R\gg 1$ so that
\begin{equation}\label{eq:n_F}
n_F(R) < R^\alpha,
\end{equation}
assume that $|z|\le \tfrac12\, R$, and set
\[
\pi_R(z) = \prod_{|\lambda|\le R} \Bigl( 1 - \frac{z}\lambda \Bigr),
\]
here and elsewhere, by $\lambda$ we denote the points of the zero set  $\Lambda_F$ of $F$ 
taken with multiplicities.
Then
\[
F(z) = \exp[ a_Rz + h_R(z)]\, \pi_R (z), 
\]
where 
\[
h_R(z) = \sum_{|\lambda| >R} \Bigl[ \log\Bigl( 1-\frac{z}{\lambda} \Bigr) + \frac{z}\lambda \Bigr] 
= - \sum_{j\ge 2} \frac{s_j (R)}j\, z^j\,, 
\]
and
\[
s_j (R) = \sum_{|\lambda| > R} \lambda^{-j}\,. 
\]
To simplify writing, we assume that $a_R\ge 0$.

\subsubsection{An outline}
By~\eqref{eq:n_F}, we have an 
a priori bound
\begin{equation}\label{eq:pi1}
\log |\pi_R(z)| \lesssim R^\alpha\log R.
\end{equation}
The chief observation is that on the positive ray we can get a better bound. Since
\[
\pi_R(z) = F(z)e^{-z} e^{-(a_R z + h_R(z)-z)},
\]
taking into account~\eqref{eq:10}, we have
\[
\log |\pi_R(r)| \le | \max_{|z|=r} (a_R z + h_R(z) ) -r | + O(1).
\]
The central part of the proof is the estimate of the RHS given in Lemma~\ref{lemma:main}.
It yields  that, for any
$\alpha~<~\alpha'~<~\tfrac12$, $\delta>0$, and $0\le r \le R^{1-\delta}$, we have
\begin{equation}\label{eq:pi2}
\log |\pi_R(r)| \le O(R^{\alpha'})\, \frac{r}{R^{1-\delta}} +O(1).
\end{equation}
Then. using the Phragm\'en--Lindel\"of idea, we consider the function
\[
v_R(z)  = \log |\pi_R(z)| - \eta r^\mu \cos\mu(\pi-|\theta|), \quad z=re^{{\rm i}\theta}, \,
a'+\delta < \mu < \tfrac12, \, \eta>0.
\]
It is subharmonic in the disk $R^{1-\delta} \mathbb D$ cut by the positive ray, and, by estimates~\eqref{eq:pi1}
and~\eqref{eq:pi2}, is bounded by 
a positive numerical constant $C$ on the boundary of the slit disk, provided that $R$ is sufficiently large. 
Hence, by the maximum principle,
$\log |\pi_R(z)| \le  \eta r^\mu \cos\mu(\pi-|\theta|) +C$ for $|z|\le R^{1-\delta}$.
Thus, $\{ \pi_R \}$ is a normal family of polynomials. Let
$\pi$ be a limiting entire function for this family. The function $\pi$ satisfies $\pi (0)=1$ and   
$\log |\pi (z)| \le  \eta r^\mu \cos\mu(\pi-|\theta|) +C$, $z\in\mathbb C$,  and vanishes on $\Lambda_F$.
Letting $\eta\to 0$, we see that $\pi$ identically equals $1$.
Hence, $\Lambda_F = \emptyset$.

\subsubsection{Preliminary estimates}
We begin with several simple estimates.

\begin{claim}\label{claim1}
For $R\to\infty$, $s_j(R) = O(R^{1-j})$, uniformly in $j\ge 2$.
\end{claim}

\noindent{\em Proof of Claim~\ref{claim1}}: Since $F$ has a finite exponential type, we have
$n_F(t) \lesssim t$, for $t\ge 1$, whence,
\[
|s_j(R)| 
\le \int_R^\infty \frac{{\rm d}n_F(t)}{t^j} \le j\, \int_R^\infty \frac{n_F(t)}{t^{j+1}}\, {\rm d}t 
\lesssim j\, \int_R^\infty \frac{{\rm d}t}{t^{j}}, {\rm d}t 
= 
\frac{j}{j-1}\, R^{1-j},
\]
proving the claim. \hfill $\Box$
 
\begin{claim}\label{claim2}
For $0\le r=|z| \le \tfrac12\, R$, we have
$ |h_R(z)| \lesssim \displaystyle \frac{r^2}R $.
\end{claim}

\noindent This claim is a straightforward consequence of the previous one. \hfill $\Box$

\begin{comment}
\medskip
Set 
\[
\pi_R^*(z) = \prod_{|\lambda|\le R} \Bigl( 1 + \frac{z}{|\lambda|} \Bigr),
\]
and observe that 
\[
 |\pi_R^*(-r)| = \min_{|z|=r} |\pi_R^*(z)| \le  
\min_{|z|=r} |\pi_R(z)|  \le \max_{|z|=r} |\pi_R(z)| \le \max_{|z|=r} |\pi_R^*(z)| = |\pi_R^*(r)|\,.
\]
\end{comment}

\begin{claim}\label{claim3}
For $0\le |z| \le \tfrac12 R$, $R\gg 1$, we have 

\smallskip\noindent {\rm (i)} 
$\log_+ |\pi_R (z)| \lesssim R^\alpha \log R$.

\smallskip\noindent {\rm (ii)} 
If, in addition, $\operatorname{dist}(z, \Lambda_F) \ge R^{-2}$,
then also $\log_- |\pi_R (z)|\lesssim R^\alpha \log R$.
\end{claim}

\noindent{\em Proof of Claim~\ref{claim3}}:
The first estimate is straightforward since each factor in the product $\pi^*_R$ is bounded by $C_FR$,
while the number of factors does not exceed $R^\alpha$.
To prove the second bound, we observe that, for $|z-\lambda|\ge R^{-2}$, we have
\[
\Bigl| 1 - \frac{z}\lambda \Bigr| = \frac{|z-\lambda|}{|\lambda|} \ge \frac1{R^2 |\lambda|}
\ge R^{-3},
\]
and again take into account the number of factors does not exceed $R^\alpha$. \hfill $\Box$

\begin{claim}\label{claim4}
For $0 \le r \le \tfrac12 R$, we have
\[
\max_{|z|=r} {\rm Re}[a_R z+h_R(z)] = r + O(R^\alpha \log R)\,.
\]
\end{claim}

\noindent{\em Proof of Claim~\ref{claim4}}:
The lower bound follows from~\eqref{eq:10} and Claim~\ref{claim3}(i):
\begin{align*}
\max_{|z|=r} {\rm Re}[a_R z+h_R(z)] &= \max_{|z|=r} \bigl[ \log |F(z)| - \log|\pi_R(z)| \bigr] \\
&\ge \log M_F(r) - \max_{|z|=r} \log |\pi_R(z)| \\
&\ge r - \tfrac12 \log r - O(1) - O(R^\alpha \log R) \\
&= r - O(R^\alpha \log R). 
\end{align*}
To get the upper bound, we set $\mathcal E_R = \bigcup_{|\lambda|\le R} D(\lambda, R^{-2})$,
and denote by $\mathcal E_R'$ the circular projection of $\mathcal E_R$ on the positive ray.
For $r\notin \mathcal E_R'$, using~\eqref{eq:10} and Claim~\ref{claim3}(ii), we get
\begin{equation}\label{eq:25}
\max_{|z|=r} {\rm Re}[a_R z+h_R(z)] \le \log M_F(r) - \min_{|z|=r} \log |\pi_R(z)|
\le r + O(R^\alpha \log R).
\end{equation}
Note that, for $R\to\infty$,
$ \sum_{|\la|\le R} R^{-2}  = n_F(R) R^{-2} = o(1)$,
and therefore, the set $\mathcal E_R'$ cannot contain intervals of length bigger than $1$, provided
that $R$ was chosen sufficiently large. Hence,~\eqref{eq:25}
holds for all $r\in [0, \tfrac12 R]$. \hfill $\Box$

\begin{claim}\label{claim5} $|a_R - 1| \ll 1$, provided that $R$ is sufficiently large.
\end{claim}

\noindent This is a straightforward consequence to Claims~\ref{claim2} and~\ref{claim4}.
\hfill $\Box$

\medskip 

It is worth mentioning that combining these preliminary estimates 
with the idea described in the outline, one gets a weaker version of Theorem~\ref{thm:main2}(ii) with 
the exponent $\tfrac12$ replaced by $\tfrac13$.

\subsubsection{A Lemma}
To proceed, we need a better estimate for the remainder in Claim~\ref{claim4}. 
\begin{itemize}
\item
In what follows, when we let $R\to\infty$, we always assume that $R$ runs to infinity through the set
of values, for which~\eqref{eq:n_F} holds.
\end{itemize}

\begin{lemma}\label{lemma:main}
Given $\alpha<\alpha'<\tfrac12$ and a sufficiently small $\delta>0$, 
for any $r~\in~[0,  \displaystyle R^{1-\delta}]$, we have
\begin{equation}\label{eq:lemma}
\bigl| \max_{|z|=r} {\rm Re}[a_R z + h_R(z)] - r \bigr| = O (R^{\alpha'}) 
\frac{r}{R^{1-\delta}}\,, \qquad R\to\infty. 
\end{equation}
\end{lemma}

\noindent {\em Proof of Lemma~\ref{lemma:main}}:
Set
\begin{align*}
S_{r, R}(\theta) &=   {\rm Re}[a_R re^{{\rm i}\theta} + h_R(re^{{\rm i}\theta})] \\
&=  a_R r\cos\theta - \sum_{j\ge 2} \frac{|s_j|}j r^j \cos(j\theta + \theta_j),
\qquad \theta_j = {\rm arg} (s_j).
\end{align*}
Note that, for $R\to\infty$,
\begin{align*}
S_{r, R}(\theta) &=  \bigl[ a_R \cos\theta + O(r/R)\bigr] r  \qquad\qquad &({\rm by\ Claim~\ref{claim2}}) \\
&= \bigl[  \cos\theta + o(1) \bigr] r  \qquad\qquad\   &({\rm by\ Claim~\ref{claim5}}).
\end{align*}
So, for sufficiently large $R$ and for $r\le R^{1-\delta}$,
the function $S_{r, R}(\theta)$ attains its maximal value near $\theta = 0$. 

Next, we note that
\begin{align*}
S^{\prime\prime}_{r, R}(\theta) &= 
- a_R r\cos\theta + \sum_{j\ge 2} j |s_j| r^j \cos(j\theta+\theta_j)
\\
&\le -r \bigl[ c_1 + O(1)\, \sum_{j\ge 1} j (r/R)^j \bigr] \qquad\qquad {\rm for\ } |\theta|\le \theta_0 < \frac{\pi}2 \\
&= -r \bigl[ c_1 + O(r/R) \bigr] <0,
\end{align*}
provided that $\theta_0$ was chosen sufficiently small. Thus, the function $S_{r, R}(\theta)$ has a unique point of maximum, which we denote by $\vartheta_R(r)$, and which is to be found from the equation $S^\prime_{r, R}(\theta) = 0$, that is,
from the equation
\begin{equation}\label{eq:theta}
a_R \sin\theta + \sum_{j\ge 2} |s_j|r^{j-1} \sin(j\theta+\theta_j) = 0.
\end{equation}
Now,
\begin{itemize}
\item
We complexify both variables $r$ and $\theta$, and consider this equation in the polydisk
$\Delta=R^{1-\delta}\,\mathbb D \times \mathbb D$.
\end{itemize}

For $R\to\infty$, the sum on the LHS of equation~\eqref{eq:theta} remains small uniformly in $r$ and $\theta$ in the polydisk $\Delta$. Indeed, by Claim~\ref{claim1},
\[
|s_j| |r|^{j-1} \lesssim \Bigl( \frac{|r|}R \Bigr)^{j-1} \le R^{(1-j)\delta},
\]
while
\[
| \sin(j\theta+\theta_j) | \le e^j\,.
\]
So the absolute value of the sum is $o(1)$, as $R\to\infty$,
while, by Claim~\ref{claim5},
$a_R |\sin\theta|$ stays bounded from below for $|\theta|=1$. 
Thus, by Rouch\'e's theorem, equation~\eqref{eq:theta} has a unique solution 
$\vartheta_R (r)$, 
which satisfies $\vartheta_R (0) = 0$. 

Furthermore, by the implicit function theorem for analytic functions, the function $\vartheta_R$ 
is analytic in the disk $R^{1-\delta}\,\mathbb D$.
Indeed, denote by  $\Phi_R(\theta, r)$ the LHS of equation~\eqref{eq:theta}.
Then 
\[
\partial_\theta \Phi_R(\theta, r) 
= a_R \cos \theta - \sum_{j\ge 2} j |s_j| r^{j-1}
\cos(j\theta + \theta_j).
\]
For $\theta\in\mathbb D$, $r\in R^{1-\delta}\,\mathbb D$, the absolute value of the RHS 
is bounded from below by \[ c_2 - O\Bigl( \frac{|r|}R \Bigr) \ge \tfrac12\, c_2 > 0, \] provided that $R$ is
sufficiently large. Then, by the residue theorem, for sufficiently small $\delta>0$, we have
\begin{equation}\label{eq:vartheta}
\vartheta_R (r) = 
\frac1{2\pi{\rm i}} \int_{|\zeta|=\delta} 
\frac{\zeta \partial_\theta \Phi_R(\zeta, r)}{\Phi_R(\zeta, r)}\, {\rm d}\zeta\,,
\end{equation}
with the RHS analytic in $r$. Note that representation~\eqref{eq:vartheta} also yields that 
{\em the function $\vartheta_R$ remains bounded  in the 
disk $R^{1-\delta}\,\mathbb D$,  uniformly in $R$ satisfying~\eqref{eq:n_F}}. 

Next, for $r\in R^{1-\delta}\,\mathbb D$, we set
\[
g_R(r) = \Bigl[  \bigl(a_R \cos\vartheta_R(r) - 1 \bigr) - 
\sum_{j\ge 2} |s_j| r^{j-1} \cos\bigl( j\vartheta_R(r) + \theta_j \bigr) \Bigr] r,
\]
and observe that, for $0<r<R^{1-\delta}$,
\[
g_R(r) = \max_{|z|=r} {\rm Re\,} \bigl[ a_R z + h_R(z) \bigr] - r.
\]
That is, we need to show that 
\begin{equation}\label{eq:g_R}
|g_R(r)| \lesssim R^{\alpha'} \cdot \frac{r}{R^{1-\delta}}\,.
\end{equation}
We will show that this holds in the whole disk $R^{1-\delta}\mathbb D$.

The function $g_R$ is analytic in $R^{1-\delta}\mathbb D$ and vanishes at the origin.
Everywhere in $R^{1-\delta}\,\mathbb D$, we have 
\[|g_R(r)| \lesssim |r| + |r|^2/R \lesssim R^{1-\delta}, \] 
while on the radius $0\le r \le R^{1-\delta}$ we have a better bound
\[ g_R(r) = r+O(R^\alpha \log R) - r = O(R^\alpha \log R). \]

Set $I=[0, 1]$ and $D^* = \mathbb D \setminus I$,
and consider the function 
\[
\omega (\zeta) = \frac2{\pi}\, {\rm arg\,} \frac{1+\sqrt{\zeta}}{1-\sqrt{\zeta}}\,, 
\quad 0 \le {\rm arg\,}\sqrt{\zeta} \le \pi,
\]
which is harmonic in $D^*$,
equals $1$ on the circular part of $\partial D^*$, except $\zeta=1$, and vanishes on $I$, also except $\zeta=1$
(that is, $\omega(\zeta)$ is a harmonic measure of the circular part of the boundary $\partial D^*$ evaluated at 
$\zeta$).
Applying the maximum principle to the bounded subharmonic function
\[
\zeta \mapsto \log |g_R(\zeta R^{1-\delta})| - \omega(\zeta)  \sup_{R^{1-\delta} \mathbb D} \log |g|
- (1-\omega(\zeta)) \sup_{R^{1-\delta} I} \log |g|, \quad \zeta\in D^*, \]
we have
\begin{align*}
|g(\zeta R^{1-\delta})| &\le \bigl( \sup_{R^{1-\delta} \mathbb D} |g| \bigr)^{\omega(\zeta)}\, 
\bigl( \sup_{R^{1-\delta} I} |g| \bigr)^{1-\omega(\zeta)} \\
&\lesssim ( R^{1-\delta} )^{\omega(\zeta)}\, (R^\alpha \log R)^{1-\omega(\zeta)}.
\end{align*}
Then, recalling that $\omega(\zeta)$ vanishes as $\zeta\to 0$, given $\varepsilon>0$, we choose 
$R$ so big that $\omega (\zeta) < \varepsilon$ for $|\zeta|\le R^{-\delta}$, and get
\[
\sup_{|r|\le R^{1-2\delta}}\, |g_R(r)| \lesssim (R^\alpha \log R)^{1-\varepsilon} \cdot ( R^{1-\delta} )^\varepsilon \lesssim R^{\alpha'}
\]
with $\alpha<\alpha'<\tfrac12$. At last, applying Schwarz' lemma to the function $w\mapsto g_R(\zeta R^{1-2\delta})$
in the unit disk, we get~\eqref{eq:g_R} (with $2\delta$ instead of $\delta$).
\hfill $\Box$

\subsubsection{Completing the proof of Theorem~\ref{thm:main2}(ii)}
Consider the product $\pi_R(z)$
in the disk $R^{1-\delta}\mathbb D$. By Claim~\ref{claim3}(i), everywhere in the disk we have
\[
\log |\pi_R(z)| \le O(R^{\alpha'}),
\]
while,  by~\eqref{eq:10} and \eqref{eq:lemma}, we have
\begin{align*}
\log |\pi_R(r)| &\le (\log M_F(r) - r) + | \max_{|z|=r}  {\rm Re}[a_R z + h_R(z)] - r | \\
&= O(R^{\alpha'}) \frac{r}{R^{1-\delta}} + C\,,
\end{align*}
where $C$ is a positive numerical constant.

Given $\eta>0$ and $\alpha' <\mu<\tfrac12$, we take $\delta$ so small that
$\mu > \alpha' +\delta$, let $D_R$ be the disk $\{|z|<R\}$ cut along the positive ray,
and consider the function
\[
v_R(z) = \log|\pi_R(z)| - \eta r^\mu\cos\mu(\pi-|\theta|),
\]
subharmonic in $D_R$ (the ``growth-killing function'' $\eta r^\mu\cos\mu(\pi-|\theta|)$ is harmonic in $D_R$). 
On the circular part of $\partial D_R$ we have
\[
v_R (R^{1-\delta}e^{{\rm i}\theta}) \le O(R^{\alpha'}) - \eta \cos \pi\mu R^{\mu (1-\delta)} < 0,
\]
provided that $R$ is big enough (we used that $\mu (1-\delta) > \mu - \delta > \alpha'$). 
On the boundary interval we have
\[
v_R(r) \le O(R^{\alpha'})\, \frac{r}{R^{1-\delta}} + C - \eta r^\mu \cos\pi\mu.
\]
For large enough $R$, we have
$ rR^{\alpha'+\delta-1}\, \ll \,  rR^{\mu-1} \ll r^\mu $.
Hence,
$v_R(r) = O(1)$ for $0\le r \le R^{1-\delta}$. Thus, by the maximum principle,
for $|z|<R^{1-\delta}$, $R\ge R_0(\alpha', \eta, \mu)$, we get $v_R(z) \le C$, that is,
$ \log|\pi_R(z)| \le \eta r^\mu \cos\mu(\pi-|\theta|) +C $. We conclude that $\{\pi_R\}$ is
a normal family of polynomials, as $R\to\infty$ along a sequence of values with estimate~\eqref{eq:n_F}.

Let $\pi$ be a limiting entire function for the family $\{\pi_R\}$.
It equals $1$ at the origin and vanishes on $\Lambda_F$.
Furthermore, $ \log|\pi (z)| \le \eta r^\mu \cos\mu(\pi-|\theta|) + C$, $z\in\mathbb C$. 
Letting $\eta\to 0$, we conclude that $\pi$ is a constant function. Since it equals $1$ at the origin, 
$\pi=1$ identically. Thus, the zero set $\Lambda_F$ is void.
\hfill $\Box$

\section{Two examples}

\subsection{Sharpness of Theorem~\ref{thm:main2}(i)}

Here, we will show that the entire function
\[
F(z) = \sum_{n\ge 0}\, (\cos\sqrt{n}+2)\,\frac{z^n}{n!}
\]
can be written in the form $F(z)=e^z G(z)$, where $G$ is an entire function of order $\tfrac12$, and therefore,
$n_F(R)=O(\sqrt R)$ as $R\to\infty$.

Let $\phi (s)$ be the Borel-Laplace transform of the entire function $\psi (w) = \cos\sqrt w +2$, i.e.,
\[
\phi (s) = \int_0^\infty (\cos\sqrt u +2) e^{-su}\, {\rm d}u.
\]
Since $\psi $ has zero exponential type, the function $\phi$ can be analytically continued to
the closed complex plane punctured at the origin,
vanishes at $\infty$, and, by the inversion formula, for each $\rho>0$,
\[
\cos\sqrt w + 2 = \frac1{2\pi {\rm i}}\, \oint _{|s|=\rho} \phi(s) e^{s w}\, {\rm d}s.
\]
Thus,
\begin{align*}
F(z) &= \sum_{n\ge 0} \frac1{2\pi {\rm i}} \oint_{|s|=\rho} \frac{e^{ns}z^n}{n!}\, \phi(s)\, {\rm d}s \\
&=\frac1{2\pi {\rm i}} \oint_{|s|=\rho} e^{z e^s}\, \phi(s)\, {\rm d}s,
\end{align*}
and therefore,
\[
G(z)=\frac1{2\pi {\rm i}} \oint_{|s|=\rho} e^{z(e^s-1)}\, \phi(s)\, {\rm d}s.
\]

To estimate the function $G$, first, we estimate the function $\phi$ and then will optimize the value of the parameter $\rho$.
We claim that
\[
|\phi (s)| = O(1) e^{C/|s|}\,.
\]
Indeed, rotating the integration line, we get
\begin{align*}
|\phi (s) | = &\Big|\, \int_0^{\infty \bar s} (\cos\sqrt w +2) e^{-ws}\, {\rm d}w\, \Big| \\
&\le \int_0^\infty e^{\sqrt t - t|s|}\, {\rm d} t.
\end{align*}
The only critical point of the function $t\mapsto \sqrt t - t|s|$ is $t^* = (4|s|^2)^{-1}$, which is a non-degenerate maximum, and $\sqrt{t^*} - t^*|s| = (4|s|)^{-1}$. Hence, by a standard Laplace asymptotic estimate,
$|\phi (s) |\le O(1)\, e^{C/|s|} $, proving the claim.

Now,
\begin{align*}
|G(z)| &\le O(1)\, \int_{-\pi}^\pi e^{|z|(e^{\rho e^{{\rm i}\theta}}+C\rho^{-1})}\, \rho {\rm d}\theta \\
&\le O(\rho) e^{C(|z|\rho + \rho^{-1})}\,.
\end{align*}
Choosing $\rho = 1/\sqrt{|z|}$, we get $|G(z)|\le O\bigl( e^{C|z|^{1/2}} \bigr)$,
as was planned. \hfill $\Box$

\subsection{A subharmonic example}

For a subharmonic function $u$ in a planar domain, we denote by $\mu_u$ its {\em Riesz measure},
that is, $\mu_u = (2\pi)^{-1}\Delta u$, where the Laplacian is understood in the sense of distributions.

\begin{proposition}\label{prop}
There exists a subharmonic function $u\colon R\mathbb D \to \mathbb R$ with
\begin{enumerate}
\item $u(r)=r \ge \max_\theta u(re^{{\rm i}\theta}) - 5$ for $0\le r\le R$,
\item $\mu_u(\mathbb D) \simeq 1$,
\item $\mu_u(R\mathbb D) \ll \sqrt{R}$ for large $R$.
\end{enumerate}
\end{proposition}
This shows that the method of proof of Theorem~\ref{thm:main2}(ii) cannot work with a 
weaker assumption $n_F(R_j) \ll \sqrt{R_j}$ for a sequence $R_j\to\infty$ because the above properties
are the only ones used in the proof. On the other hand, we  have no idea how to extend $u$ to a subharmonic
function of linear growth in the whole plane $\mathbb C$ so that this picture repeats infinitely many times. 

\medskip\noindent{\em Construction of the function $u$}: Take $\alpha > 0$ and consider the function
$f(z)=e^{-\alpha z^2}$ in the semi-disk $\mathbb D_+ = \mathbb D \cap \{{\rm Im}(z)\ge 0\}$.

\begin{claim}\label{claimA}
$|{\rm Im} f(z)| < 5 e^{2\alpha}y^2 + 2 \alpha |x| e^{-\frac12 \alpha x^2} y$.
\end{claim}
\noindent{\em Proof}: Indeed,  
$ |{\rm Im} f(z)| \le e^{\alpha(y^2-x^2)}|\sin (2\alpha xy)| \le e^{\alpha(y^2-x^2)}2\alpha |xy|$. 
For $|x|\le 2y$, the RHS does not exceed 
\[
\alpha e^{\alpha(y^2-x^2)}(x^2 + y^2) \le 
5 \alpha e^\alpha y^2 < 5 e^{2\alpha}y^2,
\] 
while,
for $|x|\ge 2y$, the RHS does not exceed $2 \alpha |x| e^{-\frac34 \alpha x^2} y$.
\hfill $\Box$

\medskip
Let $F(z)=\displaystyle \int_0^z f(t)\, {\rm d}t$.
\begin{claim}\label{claimB}
${\rm Im\,} F(x+{\rm i}y) < 10 e^{2\alpha}y^2 + e^{-2\alpha} $.
\end{claim}
\noindent{\em Proof}: We have
\[
F({\rm i}y) = {\rm i} \int_0^y e^{\alpha t^2}\, {\rm d}t.
\]
Then,
\[
{\rm Im\,} F({\rm i}y) = \int_0^y e^{\alpha t^2}\, {\rm d}t 
< \int_0^y (1+e^\alpha t^2)\, {\rm d}t < y + e^\alpha y^2\,,
\]
and therefore,
\begin{align*}
{\rm Im\,} F(x+{\rm i}y) &= {\rm Im\,} F({\rm i}y) +\int_0^x {\rm Im\,} f(t+{\rm i}y)\, {\rm d}t \\
&< y + e^\alpha y^2 + 5e^{2\alpha} y^2 + 2 \alpha y \int_0^\infty xe^{-\frac12 \alpha x^2}\, {\rm d}x 
\qquad \qquad \qquad ({\rm by\ Claim~\ref{claimA}}) \\
&=y +  e^\alpha y^2 + 5e^{2\alpha} y^2 + 2y < 3y + 6 e^{2\alpha} y^2 
< 10 e^{2\alpha}y^2 + e^{-2\alpha} 
\end{align*}
(at the very end we used that $y < e^{2\alpha} y^2 +e^{-2\alpha}$).
\hfill $\Box$

\medskip
Now, we put
\[
u(z) = 
\begin{cases}
{\rm Re\,} z + \sqrt{R}\,{\rm Im\,} F(\sqrt{z/R}), & z\in R \mathbb D_+, \\
u(\bar z),  &\bar z \in R \mathbb D_+.
\end{cases}
\]
The function $u$ is continuous in $R\mathbb D$ and harmonic in $\{z\in R\mathbb D\colon {\rm Im\,}z >0\}$
as well as in  $\{z\in R\mathbb D\colon {\rm Im\,}z < 0\}$. Moreover, $\partial_y F(x+{\rm i}0) = {\rm i} f(x)$, so
$\partial_y {\rm Im\,}F(x+{\rm i}0) = f(x) >0$, and therefore, the function $u$ is subharmonic in $R\mathbb D$.
Its Riesz measure $\mu_u$ is supported by $[0, R]$ and has the linear density $2\partial_y u(x+{\rm i}0) = 2f(x)$.

Clearly, $u(r)=r$ for $0\le r \le R$. Let $z=re^{{\rm i}\theta}$, $0\le \theta \le 2\pi$. Then
\begin{align*}
u(re^{{\rm i}\theta}) &= r\cos\theta + \sqrt{R}\, {\rm Im\,}F\Bigl( \sqrt{\frac{r}R}\, e^{{\rm i}\theta/2} \Bigr) \\
&< 
r - 2r\sin^2 \frac{\theta}2 + \sqrt{R}\, \Bigl[ 10 e^{2\alpha} \sqrt{\frac{r}R}\, \sin^2 \frac{\theta}2 + e^{-2\alpha}\Bigr]
\qquad \qquad ({\rm by\ Claim}~\ref{claimB}).
\end{align*}
Choosing $\alpha$ so that $5 e^{2\alpha} = \sqrt{R}$, we get
$ u(re^{{\rm i}\theta}) \le r + 5 $ for all  $\theta\in [0, 2\pi]$.

It remains to estimate the Riesz measure $\mu_u$.
We have
\begin{align*}
\mu_u (R \bar{\mathbb D}) = 2\int_0^R \partial_y u(x+{\rm i}0)\, {\rm d}x 
&= 2\int_0^R \sqrt{R} \frac{\rm d}{{\rm d}x} F\Bigl( \sqrt{\frac{x}R}\,\Bigr)\, {\rm d}x \\
&= 2\sqrt{R} F(1) = 2\sqrt{R} \int_0^1 e^{-\alpha x^2}\, {\rm d}x 
\le \frac{C}{\sqrt{\alpha}} \sqrt{R} \simeq \sqrt{\frac{R}{\log R}}
\end{align*}
and 
\begin{align*}
\mu_u (\bar{\mathbb D}) = 2\int_0^1 \partial_y u(x+{\rm i}0)\, {\rm d}x 
&= 2\int_0^1 \sqrt{R} \frac{\rm d}{{\rm d}x} F\Bigl( \sqrt{\frac{x}R}\, \Bigr)\, {\rm d}x \\
&= 2\sqrt{R} F\Bigl(\sqrt{\frac1R}\,  \Bigr) = 2\sqrt{R} \int_0^{1/\sqrt{R}} e^{-\alpha x^2}\, {\rm d}x
\to 2, \ {\rm as\ } R\to\infty.
\end{align*}

\bigskip
\medskip

\noindent L.H.:
School of Mathematics, Tel Aviv University, Tel Aviv, Israel
\newline {\tt vsxh2001@gmail.com}

\smallskip\noindent M.S.:
School of Mathematics, Tel Aviv University, Tel Aviv, Israel
\newline {\tt sodin@tauex.tau.ac.il}

\end{document}